\documentclass{article}
\usepackage[latin1]{inputenc}
\newtheorem{defi}{Definition}

\newtheorem{prop}[defi]{Proposition}
\newtheorem{thm}[defi]{Theorem}

\newtheorem{rmk}[defi]{Remark}

\begin{document}
\title{\bf PHH harmonic submersions are stable
}
\author{Monica Alice Aprodu}
\date{}
\maketitle

\noindent
University Dun\u area de Jos, Domneasca 47,
6200, Galati, Romania, Email: Monica.APRODU@ugal.ro
\\\\
{\bf Abstract.}
We prove that PHH harmonic submersions are 
(weakly) stable.
\\\\
A harmonic map between Riemann manifolds is called 
{\em (weakly) stable} if the
Hessian of the energy functional is (semi)
positive definite, see, for example \cite{Urk},
Chapter 5. In particular, an energy-minimizing map
is stable. 
Lichnerowich has proved in 1970 (see \cite{Li}) 
that holomorphic maps between K\"ahler manifolds 
are (weakly) stable; away from these particular
mappings, we do not dispose of many other examples
of harmonic maps which are (weakly) stable. 

In the joint paper \cite{AAB}, we have introduced 
a class of harmonic maps, defined on a Riemann manifold, 
with value in  a K\"ahler manifold, called {\em PHH 
harmonic maps} which have a behaviour somewhat similar 
to that of holomorphic maps. Holomorphic maps between
K\"ahler manifolds are typical examples of
PHH harmonic maps, but examples of different
flavour have been found in \cite{AA}.

The aim of this paper is to prove that PHH
harmonic submersions are actually (weakly) stable,
yet another property which relates maps in
this class to holomorphic maps (compare to \cite{BBdBR}).

\bigskip

We start by recalling some basic notions and
facts from \cite{AAB}, and \cite{Lou}, see also \cite{BBdBR}. 
We start with a map $\varphi: (M^m,g_M)\rightarrow
(N^{2n},J_N,g_N)$, defined on a Riemann manifold
with value in a K\"ahler manifold. For any
point $x\in M$, we denote 
$d\varphi_x^*:T_{\varphi(x)}N\rightarrow T_xM$ 
the adjoint of the tanget map
$d\varphi_x:T_xM\rightarrow T_{\varphi(x)}N$.
The map $\varphi$ is called PHWC at $x$
({\em pseudo-horizontally weakly conformal})
if and only if $d\varphi_x\circ d\varphi_x^*$
commutes with $J_{N,\varphi(x)}$.
Naturally, $\varphi$ is called PHWC
if it is PHWC at any point $x\in M$.
If $\varphi$ is PHWC at $x$, we say that
$\varphi$ is PHH at $x$ ({\em pseudo-horizontally 
homothetic}) if and only if 
$$
d\varphi_x\left((\nabla_v^Md\varphi_x^*(J_NY))_x\right)
=J_{N,\varphi(x)}d\varphi_x\left((\nabla_v^Md\varphi_x^*(Y))_x\right)
$$
for any horizontal tangent vector $v\in T_xM$, and any 
vector field $Y$, locally defined in a neighbourhood
of $\varphi(x)$, where $\nabla^M$ is the
Levi-Civita connection on $M$ (the Levi-Civita
connections on $N$ is denoted by $\nabla^N$, and
the induced connection in the bundle $\varphi^{-1}TN$
is denoted by $\widetilde{\nabla}$). The map $\varphi$
is called PHH if it is PHH at any point $x$ of $M$. 
It is easy to see that a PHWC map is PHH if and only if
$$
d\varphi (\nabla ^M_Xd\varphi ^*(J_NY))=
J_Nd\varphi (\nabla ^M_Xd\varphi ^*(Y)),
$$
for any horizontal vector field $X$ on $M$ and
any vector field $Y$ on $N$.

The two conditions PHWC and PHH can be seen
in terms of the almost complex structure on the
horizontal bundle, defined by
$J_H=d\varphi ^{-1}\circ J_N\circ d\varphi$.
$J_H$ is the restriction of an $f$-structure on $M$,
also denoted by $J_H$,
which vanishes on the vertical distribution.
The PHWC condition is the compatibility condition of
$J_H$ with the metric on $M$, while 
pseudo-horizontally homothetic condition
is equivalent to $J_H$ being parallel 
in horizontal directions. 

\bigskip

Suppose next that the map $\varphi$ is harmonic
and submersive, and $M$ is compact. In this case,
we know from Theorem 2.1 (a), Proposition 3.1 and 
Proposition 3.3 in \cite{AAB} that
the fibres of $\varphi$ are minimal submanifolds.
Recall next that the stability of harmonic maps is controlled
by a condition on the Hessian of the energy-functional:
$$
H(E)_\varphi(V,V)\geq 0,
$$
for all $V$ a section of the bundle $\varphi^{-1}TN$.
The Hessian is computed by
(see \cite{Urk}, pp. 155) 
$$
H(E)_\varphi(V,W)=\int\limits_Mg_N(J_\varphi V,W)v_M,
$$
for all $V,W$ sections in $\varphi^{-1}TN$, where
$v_M$ is the volume form on $M$, and
$J_\varphi$ is a second order selfadjoint elliptic
differential operator acting on sections of $\varphi^{-1}TN$ in
the following way. Denote by $^NR$ the curvature
tensor field on $N$, and let $\{\varepsilon_1,...,\varepsilon_m\}$ 
be an orthogonal vector frame on $M$, and $V$ be a section in
$\varphi^{-1}TN$. Then
$$
J_\varphi V:=-\sum\limits_{i=1}^m
\frac{1}{||\varepsilon_i||^2}\left(
\widetilde{\nabla}_{\varepsilon_i}
\widetilde{\nabla}_{\varepsilon_i}
-\widetilde{\nabla}_{{\nabla}_
{\varepsilon_i}\varepsilon_i}\right)V
-\sum\limits_{i=1}^m\frac{1}{||\varepsilon_i||^2} 
^NR(V,d\varphi(\varepsilon_i))d\varphi(\varepsilon_i).
$$

The second-order elliptic differential operator
$$
\bar{\Delta}_\varphi:=-\sum\limits_{i=1}^m
\frac{1}{||\varepsilon_i||^2}
\left(\widetilde{\nabla}_{\varepsilon_i}
\widetilde{\nabla}_{\varepsilon_i}
-\widetilde{\nabla}_{{\nabla}_
{\varepsilon_i}\varepsilon_i}\right)
$$
is called the {\em rough Laplacian} of $\varphi$,
cf. \cite{Urk}, pp. 155. The second sum which appears 
in the formula defining the {\em Jacobi operator}
$J_\varphi$ is denoted by ${\cal R}_\varphi$, so
$J_\varphi=\bar{\Delta}_\varphi-{\cal R}_\varphi$.
One of the useful properties of the rough Laplacian,
which will be constantly used in the sequel is
the following, cf. \cite{Urk}, pp. 156.

\begin{prop}
The rough Laplacian $\bar{\Delta}_\varphi$
satisfies
$$
\int\limits_Mg_N(\bar{\Delta}_\varphi V,W)v_M=
\int\limits_Mg_N(\widetilde{\nabla}V,\widetilde{\nabla}W)v_M=
\int\limits_Mg_N(V,\bar{\Delta}_\varphi W)v_M,
$$
where $V$ and $W$ are sections on $\varphi^{-1}TN$, and
$$
g_N(\widetilde{\nabla}V,\widetilde{\nabla}W)
=\sum\limits_{i=1}^m\frac{1}{||\varepsilon_i||^2}
g_N(\widetilde{\nabla}_{\varepsilon_i}V,
\widetilde{\nabla}_{\varepsilon_i}W).
$$
\end{prop}

After these preparations, we arrive at the
statement of the main result of this paper.

\begin{thm}
Let $(M^m,g_M)$ be a compact Riemann manifold, $(N^{2n},J_N,g_N)$
be a K\"ahler manifold, and 
$\varphi :M\rightarrow N$ be a harmonic PHH submersion.
Then $\varphi$ is (weakly) stable.
\end{thm}
{\em Proof.} As in the proof of Theorem 4.1 of \cite{AAB},
we choose a (local) frame $\{e_1,...,e_n,$
$J_Ne_1,...,J_Ne_n\}$
in $\varphi^{-1}TN$ such that the system
$\{d\varphi^*(e_1),...,d\varphi^*(e_n),$
$d\varphi^*(J_Ne_1),...,d\varphi^*(J_Ne_n)\}$
is an orthogonal frame in the horizontal
distribution. We also choose $\{u_1,...,u_s\}$
an orthonormal basis for the vertical distribution.
We denote $E_i=d\varphi^*(e_i)$, and $E^\prime_i=
d\varphi^*(J_Ne_i)$, for all $i=1,...,n$.

With this notation, we apply the same strategy
of proof as in \cite{Urk}, pp. 172, Theorem 3.2.

For $V$ a section in $\varphi^{-1}TN$, we apply 
Proposition 1, and compute:
\begin{eqnarray}
\nonumber
H(E)_\varphi (V,V)& = & \int\limits_Mg_N(\widetilde{\nabla}V,
\widetilde{\nabla}V)v_M-\int\limits_Mg_N({\cal R}_\varphi V,V)v_M.
\end{eqnarray}

By definition
\begin{eqnarray}
\nonumber
g_N(\widetilde{\nabla}V,\widetilde{\nabla}V)&=&
\sum\limits_{i=1}^n\left(\frac{1}{||E_i||^2}
g_N(\widetilde{\nabla}_{E_i}V,\widetilde{\nabla}_{E_i}V)
+\frac{1}{||E^\prime_i||^2}
g_N(\widetilde{\nabla}_{E^\prime_i}V,
\widetilde{\nabla}_{E^\prime_i}V)\right)\\
\nonumber
&&+
\sum\limits_{j=1}^sg_N(\widetilde{\nabla}_{u_j}V,
\widetilde{\nabla}_{u_j}V).
\end{eqnarray}

Analogous to the operator
used in the proof of Theorem 3.2, Chapter 5,
\cite{Urk}, we define, for any $V\in
\Gamma\left(\varphi^{-1}TN\right)$,
the operator $DV\in \Gamma\left(\varphi^{-1}TN\otimes
{\cal H}^*\right)$, where $\cal H$ is
the horizontal distribution on $M$, by
$$
DV(X):=\widetilde{\nabla}_{J_HX}V-J_N\widetilde{\nabla}_XV,
$$
for any $X$ a horizontal vector field on $M$.

Next, we compute
\begin{eqnarray}
\nonumber
g_N(DV,DV)&=&\sum\limits_{i=1}^n\left\{
\frac{1}{||E_i||^2}g_N(DV(E_i),DV(E_i))\right.\\
\nonumber
&&+ \left.\frac{1}{||E^\prime_i||^2}g_N(DV(E^\prime_i),DV(E^\prime_i))
\right\}.
\end{eqnarray}

Since $J_HE_i=E^\prime_i$, $J_HE^\prime_i=-E_i$, and
$||E_i||=||E^\prime_i||$, we obtain
\begin{eqnarray}
\nonumber
g_N(DV,DV)&=&2\sum\limits_{i=1}^n
\frac{1}{||E_i||^2}\left(g_N(\widetilde{\nabla}_{E_i}V,
\widetilde{\nabla}_{E_i}V)+
g_N(\widetilde{\nabla}_{E_i^\prime}V,
\widetilde{\nabla}_{E_i^\prime}V)\right.\\
\nonumber
&&\left.-2g_N(\widetilde{\nabla}_{E_i^\prime}V,
J_N\widetilde{\nabla}_{E_i}V)\right)
\end{eqnarray}

Therefore
\begin{eqnarray}
\nonumber
\int\limits_M\left(g_N(J_\varphi V,V)-\frac{1}{2}
g_N(DV,DV)\right)v_M&=&\sum\limits_{i=1}^n
\int\limits_M\frac{1}{||E_i||^2}
\big(2g_N(\widetilde{\nabla}_{E_i^\prime}V,
J_N\widetilde{\nabla}_{E_i}V)\\
\nonumber&&
-g_N(^NR(V,d\varphi(E_i))d\varphi(E_i),V)\\
\nonumber&&
-g_N(^NR(V,d\varphi(E_i^\prime))d\varphi(E_i^\prime),V)\big)
v_M.
\end{eqnarray}

Next, taking into account the identities
$d\varphi(E_i^\prime)=J_Nd\varphi(E_i)$,
$d\varphi(E_i)=-J_Nd\varphi(E_i^\prime)$,
and the basic properties of the curvature
tensor field $^NR$, we obtain
$$
^NR(V,d\varphi(E_i))d\varphi(E_i)
+^NR(V,d\varphi(E_i^\prime))d\varphi(E_i^\prime)=
J_N\;^NR(d\varphi(E_i),d\varphi(E_i^\prime))V,
$$
and thus
$$
\int\limits_M\left(g_N(J_\varphi V,V)-\frac{1}{2}
g_N(DV,DV)\right)v_M
$$
$$
=\sum\limits_{i=1}^n\int\limits_M\frac{1}
{||E_i||^2}\left(2g_N(\widetilde{\nabla}_{E_i^\prime}V,
J_N\widetilde{\nabla}_{E_i}V)-g_N(J_N\;^NR(d\varphi(E_i),
d\varphi(E_i^\prime))V,V)\right)v_M.
$$

We compute 
\begin{eqnarray}
\nonumber
-g_N(J_N\;^NR(d\varphi(E_i),d\varphi(E_i^\prime))V,V)
=g_N(^NR(d\varphi(E_i),d\varphi(E_i^\prime))V,J_NV)
\end{eqnarray}
\begin{eqnarray}
\nonumber
&=
g_N(\widetilde{\nabla}_{E_i}\widetilde{\nabla}_{E_i^\prime}V
-\widetilde{\nabla}_{E_i^\prime}\widetilde{\nabla}_{E_i}V
-\widetilde{\nabla}_{[E_i,E_i^\prime]}V,J_NV)
\\
\nonumber
&=
E_ig_N(\widetilde{\nabla}_{E_i^\prime}V,J_NV)
-E_i^\prime g_N(\widetilde{\nabla}_{E_i}V,J_NV)
\\
\nonumber &
-g_N(\widetilde{\nabla}_{\nabla_{E_i}{E_i^\prime}}V,J_NV)
+g_N(\widetilde{\nabla}_{\nabla_{E_i^\prime}{E_i}}V,J_NV)
\\
\nonumber &
-g_N(\widetilde{\nabla}_{E_i^\prime}V,\widetilde{\nabla}_{E_i}J_NV)
+g_N(\widetilde{\nabla}_{E_i}V,\widetilde{\nabla}_{E_i^\prime}J_NV).
\end{eqnarray}

Similarly to \cite{Urk}, pp. 180, we define a
${\cal C}^\infty$ function $\phi$ on
$M$ by the formula:
$$
\phi:=\sum\limits_{i=1}^n\frac{1}{||E_i||^2}
\big( E_ig_N(\widetilde{\nabla}_{E_i^\prime}V,J_NV)
-E_i^\prime g_N(\widetilde{\nabla}_{E_i}V,J_NV)
$$
$$
-g_N(\widetilde{\nabla}_{\nabla_{E_i}{E_i^\prime}}V,J_NV)
+g_N(\widetilde{\nabla}_{\nabla_{E_i^\prime}{E_i}}V,J_NV)\big).
$$

Since 
$$
g_N(\widetilde{\nabla}_{E_i}V,\widetilde{\nabla}_{E_i^\prime}J_NV)
=-g_N(\widetilde{\nabla}_{E_i}J_NV,\widetilde{\nabla}_{E_i^\prime}V),
$$
we have
$$
\int\limits_M\left(g_N(J_\varphi V,V)-\frac{1}{2}   
g_N(DV,DV)\right)v_M=\int\limits_M\phi v_M.
$$

The proof of the Theorem will be concluded
if we prove 
$$
\int\limits_M\phi v_M=0.
$$
For this, we use Green's formula. We choose
$X$ a horizontal vector field on $M$ defined by
the property:
$$
g_M(X,Y)=g_N(\widetilde{\nabla}_{J_HY}V,J_NV),
$$
for any vector field $Y$ on $M$, and we prove $div(X)=\phi$.
Indeed, since the fibres of $\varphi$ are minimal,
and $X$ is horizontal, it follows:
$$
div(X)=\sum\limits_{i=1}^n
\frac{1}{||E_i||^2}\left(g_M(E_i,\nabla_{E_i}X)
+g_M(E_i^\prime,\nabla_{E_i^\prime}X)\right).
$$
Next,
$$
div(X)=\sum\limits_{i=1}^n\frac{1}{||E_i||^2}
\left(E_ig_M(E_i,X)-g_M(\nabla_{E_i}E_i,X)
+E_i^\prime g_M(E_i^\prime,X)
-g_M(\nabla_{E_i^\prime}E_i^\prime,X)\right)
$$
$$
=\sum\limits_{i=1}^n\frac{1}{||E_i||^2}
\left(E_ig_N(\widetilde{\nabla}_{J_HE_i}V,J_NV)
-g_N(\widetilde{\nabla}_{J_H\nabla_{E_i}E_i}V,J_NV)\right)
$$
$$
+\sum\limits_{i=1}^n\frac{1}{||E_i||^2}\left(E_i^\prime
g_N(\widetilde{\nabla}_{J_HE_i^\prime}V,J_NV)
-g_N(\widetilde{\nabla}_{J_H\nabla_{E_i^\prime}E_i^\prime}V,J_NV)\right).
$$

By the PHH condition, we have
$J_H\nabla_{E_i^\prime}E_i^\prime=-(\nabla_{E_i^\prime}E_i)^h$,
and $J_H\nabla_{E_i}E_i=(\nabla_{E_i}E_i^\prime)^h$, 
where by $(.)^h$ we denoted the horizontal component
of $(.)$, so,
$$
div(X)=\phi.
$$

We proved 
$$
\int\limits_Mg_N(J_\varphi V,V)v_M=\frac{1}{2}g_N(DV,DV)v_M\geq 0.
$$

\begin{rmk} Our result improves
the main result of \cite{Mo}, provided
that the source manifold is compact
(condition which is not needed in \cite{Mo}).
\end{rmk}

\noindent
{\em Acknowledgements.} This work was financed
by a NATO fellowship. The author is grateful
to the Fourier Institute in Grenoble for hospitality,
and to J. C. Wood for some useful remarks on
an early version of the manuscript.

\end{document}